\documentclass[12pt, letterpaper]{article}
\usepackage[OT4]{fontenc}
\usepackage{graphics,graphicx}
\usepackage{amsmath}
\usepackage{amssymb,amsfonts}
\usepackage {amsthm}
\usepackage{xcolor}

\usepackage{enumerate}
\usepackage{hyperref}
\usepackage[T2A]{fontenc}
\usepackage[cp1251]{inputenc}
\usepackage[all]{xy}
\usepackage[linesnumbered,boxed]{algorithm2e}

\newtheorem{theorem}{Theorem}
\newtheorem{proposition}{Proposition}

\newtheorem{corollary}{Corollary}
\newtheorem*{proposition*}{Proposition}

\newcommand\myfootnote[1]{
\renewcommand{\thefootnote}{}
\footnotetext{#1}
\def\thefootnote{\@arabic\c@footnote}
}

\renewcommand{\subsection}{\@startsection{subsection}{2}{0mm}{-\baselineskip}{-5pt}{\it \bf}}
\makeatother

\title{The unitary Cayley graph of upper triangular matrix rings}
\author{Waldemar Ho{\l}ubowski$^{1}$, Sergiy Kozerenko$^{2,3,4}$, \\ Bogdana Oliynyk$^{1,2}$, Viktoriia Solomko$^{1,2}$}

\begin{document}

\maketitle
	
\small \noindent 
$^{1}$Faculty of Applied Mathematics, Silesian University of Technology, Poland \\  
$^{2}$Department of Mathematics, National University of Kyiv-Mohyla Academy, Ukraine \\
$^{3}$Institute of Mathematics of NAS of Ukraine, Ukraine \\
$^{4}$Kyiv School of Economics, Ukraine\\

w.holubowski@polsl.pl, kozerenkosergiy@ukr.net, bogdana.oliynyk@polsl.pl, \\ vsolomko@polsl.pl


{\bf 2020}{ {\bf Mathematics Subject Classification:} 05C25, 05E16,  15B33, 16S50}

%05E16(2020–now)Combinatorial aspects of groups and algebras
% 15B33(2010–now)Matrices over special rings (quaternions, finite fields, etc.)
% 05C25(1973–now)Graphs and abstract algebra (groups, rings, fields, etc.)

% 16S50 Matrix rings

\begin{abstract}
	 The unitary Cayley graph $C_R$ of a finite unital ring  $R$ is the simple graph with vertex set $R$ in which  two elements $x$ and $y$ are connected by an edge if and only if $x-y$ is a unit of $R$. We characterize the unitary Cayley graph $C_{T_n (\mathbb{F})}$ of the ring of all upper triangular matrices  $T_n(\mathbb{F})$ over a finite field $\mathbb{F}$. We show that $C_{T_n (\mathbb{F})}$ is isomorphic to the semistrong product of the complete graph $K_m$ and the antipodal graph of the Hamming graph $A(H(n,p^k))$, where  $m=p^{\frac{kn(n-1)}{2}}$ and $|\mathbb{F}|=p^k$. In particular, if  $|\mathbb{F}|=2$, then the graph  $C_{T_n (\mathbb{F})}$ has $2^{n-1}$ connected components, each component is isomorphic to the complete bipartite graph $K_{m,m}$, where $m=2^{\frac{n(n-1)}{2}}$. 
	 We also compute the diameter, triameter, and clique number of the graph $C_{T_n (\mathbb{F})}$. 
\end{abstract}

%optional
%\begin{center}
%Joint work with .
%\end{center}

%\vspace{1cm}

%     54B35       Spectra ( General topology)
  %  54E35       Metric spaces, metrizability
%54E40       Special maps on metric spaces

\section*{Introduction}

In algebraic graph theory, various graphs associated with algebraic structures play a  prominent role. A notable example is the Cayley graph, denoted as $Cay(G, S)$, of a group $G$ related to some subset $S$ of $G$. The graph $Cay(G, S)$ is connected if and only if $S$ generates $G$. Furthermore, $Cay(G,S)$ is regular, and its properties provide significant insights into the structure of finitely generated groups. We refer to books~\cite{Godsil} on algebraic graph theory and~\cite{Harpe, Loh} on geometric group theory for the application of Cayley graphs in groups. 

Let $R$ be a finite ring with identity. The {\it unitary Cayley graph} of a ring $R$, denoted by $C_R$, is a simple graph with vertex set $R$, and two elements $x$, $y$ are adjacent if and only if $x-y$ is a unit of  $R$. Essentially, it is just a Cayley graph of an additive group of ring $R$ related to the set of units of $R$.

The unitary Cayley graphs have been studied intensively for the ring $Z_n$ (see, for example~\cite{Berrizbeitia, Dejter, Fuchs, Ili, Klotz07}). In particular, Dejter and Giudici~\cite{Dejter} proved that  if $n$ is a prime then the unitary Cayley graph $C_{Z_n}$ is isomorphic to a complete graph with $n$ vertices, and  $C_{Z_n}$ is a bipartite graph when $n$ is even, moreover, if $n=2^s$ then $C_{Z_n}$ is isomorphic to a complete bipartite graph $K_{2^{s-1},2^{s-1}}$. In~\cite{Akhtar} the diameter, girth, eigenvalues, vertex and edge connectivity, and vertex and edge chromatic number were described precisely for a finite commutative ring $R$. 

Let $\mathbb{F}$ be a finite field,  $n$ be a positive integer, $n \ge 2$.  In their works~\cite{Kiani1,Kiani} Kiani, Mollahajiaghaei, and Aghaei described the chromatic number, clique number, independence number, diameter  of unitary Cayley graph $C_{M_n(\mathbb{F})}$ of matrix ring $M_n(\mathbb{F})$ over a finite field $\mathbb{F}$ and proved that the graph $C_{M_n(\mathbb{F})}$ is regular and strongly regular for $n = 2$.
In a related study~\cite{Rattanakangwanwong2} Rattanakangwanwong and Meemark computed the clique number, chromatic number, and independence number of the subgraph of unitary Cayley graph of a matrix algebra over a finite field induced by the set of idempotent matrices. For some other recent papers on unitary Cayley graphs of matrix ring $M_n(\mathbb{F})$, we refer the reader to~\cite{ChenHuang, Huang, Rattanakangwanwong}.
 
In this paper, we characterize the properties of the unitary Cayley graph $C_{T_n (\mathbb{F})}$ of  the ring of all upper triangular $n \times n$ matrices $T_n(\mathbb{F})$ over the field $\mathbb{F}$. Denote by $(T_n(\mathbb{F}))^{*}$ the set of all invertible upper triangular $n \times n$ matrices. The vertices of the graph $C_{T_n (\mathbb{F})}$ are upper triangular $n \times n$ matrices $T_n(\mathbb{F})$ over  $\mathbb{F}$ and two vertices $a$ and $b$ are adjacent if and only if $a-b$ is an element of $(T_n(\mathbb{F}))^{*}$.

We will demonstrate that when $|\mathbb{F}|=2$, the unitary Cayley graph $C_{T_n (\mathbb{F})}$ of the ring of all upper triangular matrices  $T_n(\mathbb{F})$ over  $\mathbb{F}$ exhibits a fairly simple structure. 

\begin{theorem}\label{Theorem 1}
	If  $|\mathbb{F}|=2$, then the graph  $C_{T_n (\mathbb{F})}$ has $2^{n-1}$ connected components, and its each component is isomorphic to the complete bipartite graph $K_{m,m}$, where $m=2^{\frac{n(n-1)}{2}}$.				  
\end{theorem}

In the case of $|\mathbb{F}|>2$, the unitary Cayley graph $C_{T_n (\mathbb{F})}$ is connected.  Let $S$ be a set of $p$ elements, $p$ and $l$ be positive integers. The {\it Hamming graph} $H(l,p)$ is defined on the set of ordered $l$-tuples of elements of $S$ with two vertices being adjacent if they differ in precisely one coordinate. Hamming graphs find active applications in information theory and computer science. Specifically, we will consider Hamming  graph $H(n,p^k)$ for $S=\mathbb{F}$,  $|\mathbb{F}|=p^k$.

The antipodal graph of a graph $G$, denoted by $A(G)$, has the same vertex set as $G$ with two vertices $u$ and $v$ being adjacent if the distance between $u$ and $v$ is equal to the diameter of $G$ (see \cite{Aravamudhan, Johns}). 
The antipodal graph $A(H(n,p^k))$ of the Hamming graph $H(n,p^k)$  is defined on the set of all $n$-tuples of elements from the field $\mathbb{F}$, and two $n$-tuples $u=(u_1, \ldots, u_n)$ and $v=(v_1, \ldots, v_n )$ are adjacent provided $u_i \ne v_i$ for every $1 \le i \le n$.

For the case of $|\mathbb{F}|>2$, we prove the following  result that describes  the structure of graph $C_{T_n (\mathbb{F})}$ as a semistrong product of two graphs.

\begin{theorem}\label{Theorem3}
	Let $|\mathbb{F}|>2$, $|\mathbb{F}|=p^k$ for a prime $p$ and a positive integer $k$.  Then
	\begin{equation}
		C_{T_n(\mathbb{F})} \simeq K_m \bullet A(H(n,p^k)),
	\end{equation}
	where $m=p^{\frac{kn(n-1)}{2}}$.
\end{theorem}

 Additionally, we also characterize the diameter, triameter, and clique number of the graph $C_{T_n (\mathbb{F})}$. 

\section{Preliminaries}

In this paper, we consider only finite simple graphs.  
Let $G=(V(G), E(G))$ be a finite connected simple graph. Define a metric $d_{G}$ on the set of vertices $V(G)$  as follows: for any $u, v \in V(G)$ the distance $d_G(u,v)$ is defined as the length of the shortest path between $u$ and $v$.

The \emph{diameter} of a connected graph $G$ is the value 
\[diam(G)=\max\{d_{G}(u,v):u,v\in V(G)\}.\] 

For every triplet of vertices $u,v,w\in V(G)$, we define 
\[d_{G}(u,v,w)=d_{G}(u,v)+d_{G}(u,w)+d_{G}(v,w).\] 

The \emph{triameter} of a connected graph $G$ is defined as the value 
\[triam(G)=\max\{d_{G}(u,v,w):u,v,w,\in V(G)\}.\]

The triplet of vertices $u,v,w\in V(G)$ is called \emph{triametral} if $d_{G}(u,v,w)=tr(G)$. The main motivation for studying $tr(G)$ comes from its appearance in lower bounds on the radio $k$-chromatic number of a graph and the total domination number of a connected graph (\cite{Das:21, HakKozOl, KoPan:15, SahPan:15}).

A {\it clique} is a subgraph of a graph $G$ that is  isomorphic to a complete graph. The {\it clique number} of $G$ is the size of the largest clique in $G$, denoted by $\omega(G)$.

Now, consider two finite simple graphs, $G$ and $H$. The semistrong product $G \bullet H$  (see \cite{Semistrong}) is the graph with vertex set $V(G \bullet H)=V(G) \times V(H)$ and edge set 
\begin{multline*} E(G \bullet H)=$ $=\{(u_1,v_1)(u_2,v_2)|u_1u_2 \in E(G)  \text{ and }  v_1v_2 \in E(H)  \text{ or } \\ u_1=u_2  \text{ and }  v_1v_2 \in E(H) \}. 
\end{multline*}	 

\section{The properties of the unitary Cayley graph of $T_n(\mathbb{F})$}

For any unital finite ring $R$ the  unitary Cayley graph $C_R$ is $|R^{*}|$-regular \cite{Klotz07}.  As $|(T_n(\mathbb{F}))^{*}|= (|\mathbb{F}|-1)^n \cdot |\mathbb{F}|^{\frac{n^2-n}{2}}$, the next proposition becomes evident.
\begin{proposition}
	\label{Prop0}
  The graph $C_{T_n (\mathbb{F})}$ is $(|\mathbb{F}|-1)^n |\mathbb{F} |^{\frac{n^2-n}{2}}$-regular. 
\end{proposition}

From the structure of upper triangular matrices, we deduce a simple condition for the existence of an edge in the graph $C_{T_n (\mathbb{F})}$.

\begin{proposition}\label{Prop1}
	 Two matrices $a =\{a_{ij}\}_{1 \le i,j \le n}$ and $b =\{b_{ij}\}_{1 \le i,j \le n}$ are connected by edge in $C_{T_n (\mathbb{F})}$ if and only if $a_{ii}-b_{ii} \ne 0$ for any $1 \le i,j \le n$.  
\end{proposition}
\begin{proof}
	Indeed, $a-b$ is invertible if and only if $det(a-b) \ne 0$. As $a$ and $b$ are upper triangular,  
	\[det(a-b)=\prod_{i=1}^{n} (a_{ii}-b_{ii}).\]
	This completes the proof of the proposition.  
\end{proof} 

\begin{theorem}\label{Theorem2}
	If $|\mathbb{F}| \ne 2$, then the graph $C_{T_n (\mathbb{F})}$ is connected.
\end{theorem}

First, we prove Theorem~\ref{Theorem 1}.

\begin{proof}[Proof of Theorem~\ref{Theorem 1}]

Let $|\mathbb{F}|=2$ and let $a =\{a_{ij}\}_{1 \le i,j \le n}$ be a matrix from $T_n(\mathbb{F})$. For the vector $(a_{11}, a_{22}, \ldots, a_{nn})$  there exists only one vector  $(b_{11}, b_{22}, \ldots, b_{nn}) $ such that  
\[a_{ii}-b_{ii}=1 \mod 2, \qquad \text{fo any $1\le i \le n$} .\] 
Proposition~\ref{Prop1}	 implies, that the matrix  $a$ is connected by an edge with any matrix  $b =\{b_{ij}\}_{1 \le i,j \le n}$ from $T_n(\mathbb{F})$ with the elements $(b_{11}, b_{22}, \ldots, b_{nn})$ on the main diagonal. It is true for all matrices with  the elements $(a_{11}, a_{22}, \ldots, a_{nn})$ on the main diagonal.  However, any matrix with the elements $(b_{11}, b_{22}, \ldots, b_{nn})$ on the main diagonal is connected by an edge with any matrix with the elements $(a_{11}, a_{22}, \ldots, a_{nn})$ on the main diagonal and is not connected with any another matrix from $T_n(\mathbb{F})$.

Hence, all matrices with the elements $(a_{11}, a_{22}, \ldots, a_{nn})$ or $(b_{11}, b_{22}, \ldots, b_{nn})$ on the main diagonal form the connected component of the graph $C_{T_n (\mathbb{F})}$. Since  we choose matrix $a$ arbitrarily, any connected component of graph $C_{T_n (\mathbb{F})}$ is defined by two vectors $(c_{11}, c_{22}, \ldots, c_{nn})$ and $(d_{11}, d_{22}, \ldots, d_{nn})$, such that 
\[c_{ii}-d_{ii}=1 \mod 2, \text{ for any } 1\le i \le n.\] 
Therefore, the graph $C_{T_n (\mathbb{F})}$ has $2^{n-1}$ connected components and any component is isomorphic to $K_{m,m}$ where $m$ equals the number of  upper triangular $n \times n$ matrices over field $\mathbb{F}_2$ with elements $(c_{11}, c_{22}, \ldots, c_{nn})$ on the main diagonal, i.e. 
\[m=2^{\frac{n(n-1)}{2}}.\]
\end{proof} 

\begin{proof}[Proof of Theorem~\ref{Theorem2}]	
Let now $|\mathbb{F}|>2$. Assume, that  matrices $a =\{a_{ij}\}_{1 \le i,j \le n}$ and $b =\{b_{ij}\}_{1 \le i,j \le n}$ are not connected by an edge. Then we can construct a matrix   $c =\{c_{ij}\}_{1 \le i,j \le n}$ in the following  way: \begin{itemize}
	\item for all $i$, $ 1 \le i \le n$, we choose $c_{ii}$ such that  $c_{ii} \ne a_{ii}$,  $c_{ii} \ne b_{ii}$ (we can choose $c_{ii}$  because  $|F|>2$);
	
	\item 
	for all $ 1 \ge i, j \ge n $, $i \ne j$ we set $c_{ij} =0$. 
\end{itemize}

Proposition~\ref{Prop1}	 implies that the matrix $c$ is connected by an edge with the matrices $a$ and $b$. So, for any matrices  $a$ and $b$ that are not connected by an edge, there exists a matrix $c$ which is connected by an edge with both $a$ and $b$ simultaneously. Therefore, the graph $C_{T_n (\mathbb{F})}$ is connected. 
\end{proof}	

This theorem directly implies the next proposition. 

\begin{corollary}
	\label{Cor1} 
	If $|\mathbb{F}|>2$, then 
	\[diam (C_{T_n (\mathbb{F})})=2.\] 
\end{corollary}

\begin{proof}
The proof of this corollary directly follows from the proof of Theorem~\ref{Theorem2}.  
\end{proof} 	

The next proposition describes the triameter of the graph $C_{T_n (\mathbb{F})}$.
\begin{proposition} 
	If $|\mathbb{F}|>2$, then 
	\[triam (C_{T_n(\mathbb{F})})=6.\] 
\end{proposition}
\begin{proof}
It is clear that  $triam (C_{T_n(\mathbb{F})}) \leq 6$, because $triam(G)=\max\{d_{G}(u,v,w):u,v,w,\in V(G)\}$ and $d_{G}(u,v,w)=d_{G}(u,v)+d_{G}(u,w)+d_{G}(v,w) \leq  3 \cdot diam(G)$.
Since $|F| > 2$, there exist three pairwise different elements $a,b, c \in F$.  We consider three  diagonal matrices $d_1= diag(a,a, a,  \ldots, a)$, $d_2=diag(a, b, b,  \ldots, b)$ and 
$d_3=diag(a, c, c, \ldots, c)$. We have   $d_{C_{T_n(\mathbb{F})}}(d_1,d_2,d_3)=2+2+2=6$. This completes the proof of the proposition.

\end{proof}

The next proposition characterizes the  clique  number of the graph $C_{T_n (\mathbb{F})}$.
\begin{proposition} 
	 \[\omega(C_{T_n(\mathbb{F})})=|\mathbb{F}|.\]  
\end{proposition}

\begin{proof}
  Let $S$ be a  maximal clique of a graph $G_{T_n}$. Proposition~\ref{Prop1} implies, that 
  \[a_{ii} \ne b_{ii}\] 
  for any two different elements    $a =\{a_{ij}\}_{1 \le i,j \le n} $ and $b =\{b_{ij}\}_{1 \le i,j \le n}$ from $S$ and for any $i$,  $1\le i \le n$. But for any $i$ we can choose  only $|\mathbb{F}|$ different elements $a_{ii}$. So, $|S| = |\mathbb{F}|$.     
 \end{proof} 	

\section{Connection with Hamming graph}

\begin{proof}[Proof of Theorem~\ref{Theorem3}]
	
We have to	show that the unitary Cayley graph $C_{T_n (\mathbb{F})}$ of the ring of all upper triangular matrices  $T_n(\mathbb{F})$ over  $\mathbb{F}$	is isomorphic to the semistrong product of the complete graph $K_m$ and the antipodal graph of the Hamming graph $A(H(n,p^k))$, where  $m=p^{\frac{kn(n-1)}{2}}.$
	
First, we define a complete graph on the set of all strictly upper triangular matrices. This graph is isomorphic to the graph $K_m$, where $m=p^{\frac{kn(n-1)}{2}}$. 
	
Now, determine a bijection $\varphi$ from the set of vertices of the graph $C_{T_n (\mathbb{F})}$ to the Cartesian product of the set of vertices of $K_m$ and the set of vertices of the graph $A(H(n,p^k))$.  Let  $a =\{a_{ij}\}_{1 \le i,j \le n}$ be a matrix from $T_n(\mathbb{F})$. Define
	\[\varphi(a) = (\hat{a} , (a_{11}, \ldots, a_{nn})),\]
where 
\[\hat{a}  =\begin{pmatrix}
	0 &   a_{12} & \ldots  & a_{1n} \\
	0 &  0 & \ldots  & a_{2n} \\
	\vdots& \vdots&  \vdots & \vdots \\
	0 &   0 & \ldots  & a_{(n-1),n} \\
	0 &  0 & \ldots  & 0 
\end{pmatrix}.\]
   
We would like to show that two matrices $a =\{a_{ij}\}_{1 \le i,j \le n}$ and $b =\{b_{ij}\}_{1 \le i,j \le n}$ are connected by an edge in the graph $C_{T_n (\mathbb{F})}$ if and only if they are connected by edge in the semistrong product of graphs $K_m$ and $ A(H(n,p^k))$. Assume, that there exist the edge $ab$ in the graph $C_{T_n (\mathbb{F})}$. Proposition~\ref{Prop1} implies that $a_{ii}\ne b_{ii}$ for any $i$, $1\le i \le n$. Then the vectors $(a_{11}, \ldots, a_{nn})$ and $(b_{11}, \ldots, b_{nn})$ are connected by an edge in  $A(H(n,p^k))$.  If the matrices
 \[\hat{a}  =\begin{pmatrix}
	0 &   a_{12} & \ldots  & a_{1n} \\
	0 &  0 & \ldots  & a_{2n} \\
	\vdots& \vdots&  \vdots & \vdots \\
	0 &   0 & \ldots  & a_{(n-1),n} \\
	0 &  0 & \ldots  & 0 
\end{pmatrix} \quad \text{ and } \quad
\hat{b}  =\begin{pmatrix}
	0 &   b_{12} & \ldots  & b_{1n} \\
	0 &  0 & \ldots  & b_{2n} \\
	\vdots& \vdots&  \vdots & \vdots \\
	0 &   0 & \ldots  & b_{(n-1),n} \\
	0 &  0 & \ldots  & 0 
\end{pmatrix}\] 
 are different, then they are connected by an edge in $K_m$. Therefore, from the definition of semistrong product of graphs follows, that there exists the edge  $ab$ in the graph $K_m \bullet A(H(n,p^k))$. 
 
 Assume there is no edge $ab$ in $C_{T_n (\mathbb{F})}$. Proposition~\ref{Prop1} implies that there exist  $i$, $1\le i \le n$, such that $a_{ii}= b_{ii}$. Consequently, the vectors $(a_{11}, \ldots, a_{nn})$ and $(b_{11}, \ldots, b_{nn})$ are not connected by edge in $A(H(n,p^k))$. Thus, there is no edge $ab$ in $K_m \bullet A(H(n,p^k))$. 
   
\end{proof}

We can describe the connection between  $C_{T_n (\mathbb{F})}$ and   $ A(H(n,p^k))$ in another way. Define on $T_n(\mathbb{F})$ equivalence relation $\equiv$: 

{\it a matrix $a =\{a_{ij}\}_{1 \le i,j \le n}$ is equivalent to a matrix $b =\{b_{ij}\}_{1 \le i,j \le n}$ if and only if   $a_{ii}= b_{ii}$ for all  $i$, $1\le i \le n$.}

That is, each equivalence class is determined by the main diagonal of  matrices. 

Let $\hat{C}_{T_n}$  be a graph induced by a graph  $C_{T_n (\mathbb{F})}$ on the set  $T_n(\mathbb{F})|_{\equiv}$. Then we have the next description of  $\hat{C}_{T_n}$.

\begin{proposition} 
	 $\hat{C}_{T_n}\simeq A(H(n,p^k))$, for all $n \ge 2$.   
\end{proposition}

The proof of this proposition  directly follows from the Theorem~\ref{Theorem3}. 

\section{Acknowledgments}

Sergiy Kozerenko is deeply thankful to the Armed Forces of Ukraine for ensuring the safety of Kyiv during the work on this paper. Viktoriia Solomko was supported by a grant from the Institute of International Education’s Scholar Rescue Fund.

\end{document}